\documentclass{amsart}

\numberwithin{equation}{section}
\newtheorem{theo}{Theorem}

\newtheorem{lem}[equation]{Lemma}
\newtheorem{cor}[theo]{Corollary}
\newtheorem{smallcor}[equation]{Corollary}

\theoremstyle{remark}
\newtheorem{rem}[equation]{Remark}
\newtheorem{rems}[equation]{Remarks} 
\newtheorem{conj}[theo]{Conjecture}
\newtheorem{conj*}{Conjecture}

\theoremstyle{definition}
\newtheorem{define}{Definition}

\newcommand\R{{\mathbf{R}}}
\newcommand\Z{{\mathbf{Z}}}  
\renewcommand\d{\partial}

\newcommand\vol{\operatorname{vol}}

\newcommand\del{\nabla}
\newcommand\e{\epsilon}

\renewcommand\to{\rightarrow}

\title[Harmonic maps]{Deforming a map into a harmonic map}
\author{Deane Yang}
\address{Department of Applied Mathematics and Physics\\
Polytechnic University\\ Six Metrotech Center\\
Brooklyn NY 11201}
\email{yang@math.poly.edu}
\date{\today}

\thanks{I was partially supported by National Science Foundation grant DMS-9200576.
 Some of the work in this paper was done at l'Institut des
Hautes Etudes Scientifiques. I would like thank the
 staff  and  the  director,  Jean--Pierre Bourguignon, for their support and hospitality.
I would also like to thank Stephen Semmes, Curt McMullen, and Michael Wolf for
helpful discussions. I am grateful to Peter Li for his comments on an 
earlier version of this paper and references to
related results.}

\begin{document}
\maketitle

\newcommand\phat{\hat{p}}
\newcommand\Bhat{\widehat{B}}
\newcommand\phihat{\hat{\phi}}
\renewcommand\H{\mathbf{H}}
\newcommand\qhat{\hat{q}}

\section{Introduction}

Let
$X$ be
 a complete noncompact Riemannian manifold with Ricci curvature and 
 Sobolev radius (see \S\ref{basic} for the definition)
 bounded from below and $Y$ a complete
Riemannian manifold with nonpositive sectional curvature.
We shall study some situations where a smooth map $f: X \to Y$ can be
deformed continuously into a harmonic map, using a naturally defined
flow. The flow used here is not the usual harmonic heat flow, as
introduced by Eells--Sampson. We use, instead, a flow introduced by
J.P. Anderson \cite{and91}.

Except for some classical results on linear elliptic partial differential 
equations, this paper is self--contained and 
provides a straightforward proof for a wide range
of existence and uniqueness theorems for harmonic maps.  
In particular, we obtain as a corollary a recent result of Hardt--Wolf \cite{HarWol96} 
on the existence of harmonic quasiisometries of the hyperbolic plane.

\section{Statement of theorems}
Throughout this paper we shall assume that $X$ is an $n$--dimensional complete 
noncompact Riemannian manifold with  Ricci curvature bounded from below by 
$-(n-1)\kappa^{2} \le 0$ and Sobolev radius bounded from below 
by $2\rho > 0$ (see \S\ref{basic} for the definition of Sobolev 
radius) and that $Y$ is an $m$--dimensional
complete Riemannian manifold with nonpositive sectional curvature.
Also, see \S\ref{defines} for
other definitions.

\begin{define}
$$
\lambda_0(X) = \inf_{f\in C^\infty_0(X)} \frac{\int_X |\del
f|^2}{\int f^2}
$$
\end{define}
Note that $\lambda_{0}(X) > 0$ only if $X$ is noncompact and has 
infinite volume.

\begin{theo} \label{unique-2}Assume $\lambda_0(X) > 0$.  
Given a harmonic map $u: X \to Y$, there is no other
harmonic map within finite $L^p$ distance of $u$, for any $1 < p <
\infty$.
\end{theo}

\begin{theo}\label{unique-1} Assume that the sectional curvature of 
$Y$ is bounded from above by $-K^{2} < 0$.
Given a strictly nondegenerate harmonic map $u: X \to Y$,
 there is no other homotopic harmonic map that is bounded distance away from $u$.
\end{theo}

\begin{theo} \label{exist-2} Assume $\lambda_0(X) > 0$.
Given $p \ge 2$ and a smooth map
$u_0: X \to Y$ satisfying
\begin{align*}
\|\Delta u_{0}\|_p &< \infty
\end{align*}
there exists a homotopic
harmonic map $u: X \to Y$ within finite  $L^p$ distance
of $u_0$.
\end{theo}

\begin{theo}\label{exist-1} 
Given $n, m, \rho, \kappa, K, \tau, C > 0$,
 there exists $\e > 0$ such that given $X$ satisfying the assumptions 
 above, $Y$ with sectional curvature
 bounded from above by $-K^{2}$, and any smooth map $u_0: X \to Y$
satisfying
\begin{align*}
\tau[u_0] &\ge \tau\\
\|\d u_{0}\|_{\infty} &< C\\
\|\Delta u_0\|_{\infty} & \le \e
\end{align*}
there exists a homotopic harmonic map $u: X \to Y$
that is a bounded distance away from $u_0$.
\end{theo}

\begin{cor} \label{exist-1a}
Given $K, C > 0$ and $\sigma > 1$, there exists $\e > 0$ such that
given any complete Riemannian manifold $(Y,g)$ with
sectional curvature bounded from above by $-K^{2}$and
any Riemannian metric $h$ on $Y$ satisfying
\begin{align*}
\sigma^{-2}g \le h  &\le \sigma^2g\\
\|\del_{g}h\|_{\infty} &\le C\\
\left\|g^{ij}\left(\del_ih_{jk} -
\frac{1}{2}\del_kh_{ij}\right)\right\|_\infty &\le \e
\end{align*}
there exists a harmonic quasiisometry $u: (X,h) \to (X,g)$ homotopic
to and a finite distance from the identity map.
\end{cor}

Let $\H^{n}$ denote hyperbolic $n$--space. By Theorem~\ref{exist-1} and Corollary~\ref{qiext},
\begin{cor}\label{exist-1b} Let $n \ge 2$. Given a harmonic 
quasiisometry $u_{0}: \H^{n} \to \H^{n}$ with boundary map 
$\hat{u}_{0}$, there exists $\delta > 0$ 
such that given any $(1+\delta)$--quasisymmetric map $\psi: 
S^{n-1} \to S^{n-1}$, the boundary map $\psi\circ\hat{u}_{0}: 
S^{n-1} \to S^{n-1}$ extends to a harmonic quasiisometry $u: \H^{n} 
\to \H^{n}$.
\end{cor}

\begin{cor}\label{exist-1c} Let $n \ge 2$. There exists $\delta > 0$, such that any
$(1+\delta)$--quasisymmetric map $\hat{u}: S^{n-1} \to S^{n-1}$
extends to a harmonic quasiisometry $u: \H^n \to \H^n$.
\end{cor}

\begin{rems}
Theorems  \ref{unique-2} and \ref{exist-2} were first proved by Ding--Wang 
\cite{DinWan91}. The proof presented here is very similar to 
theirs; the only difference is that by using the flow introduced 
here, we avoid relying on any past results on harmonic maps. Only the
existence and uniqueness of a solution to the Dirichlet problem for 
a linear elliptic system of partial differential equations is needed 
here.

Corollary \ref{exist-1a} generalizes a result of J. P. Anderson 
\cite{and91}.

Corollaries \ref{exist-1b} and \ref{exist-1c} were first proved
by Earle--Fowler in dimension 2 
and in all dimensions 
by Hardt--Wolf \cite{HarWol96} in 1993.
Another proof appears in a recent paper of Tam--Wan,  \cite{TamWan96}.
\end{rems}

\begin{rem}Douady and Earle \cite{DouEar86} found a conformally natural means of 
extending a suitable map between the boundaries of two balls to a map between 
the balls. Although they consider only the case where the domain and 
range have equal dimension and the boundary map is a homeomorphism, 
their construction works without change even if the dimensions are not 
equal and for any measurable map of the boundaries where the inverse 
image of any set with measure zero has measure zero. 
This extension is a nonlinear
analogue of the Poisson kernel for the hyperbolic Laplacian. Following
this analogy, the space of all boundary maps that extend to maps
between hyperbolic spaces with $L^p$ bounded Laplacian is a nonlinear
analogue of a Besov space. In particular, it can be shown that any
$C^{1,\alpha}$ map with everywhere nonvanishing differential
of the boundaries extends, via the Douady--Earle
extension, to a map of hyperbolic spaces with $L^p$ bounded Laplacian.
This observation and Corollary~\ref{exist-2} imply 
an existence theorem similar to those of Li--Tam
\cite{LiTam93, LiTam93A, LiTamWan95} for the Dirichlet problem at
infinity for harmonic maps between hyperbolic spaces. A detailed
discussion of this will appear elsewhere.
\end{rem}

It appears that the condition of bounded energy density is a natural
one to impose. A reasonable question to ask whether,
given Hadamard manifolds $X$ and $Y$ and a map between the respective 
boundaries at infinity,
there exists an extension of the map to a harmonic map of the interiors
with bounded energy density. It is easy, however, to give examples 
where the answer is no. For example, the map that sends the entire 
boundary of $X$ to a single point of the boundary of $Y$ has no 
harmonic extension. A more realistic conjecture appears to be the 
following:

\begin{conj} Given Hadamard manifolds $X$ and $Y$, an admissible map between the 
respective boundaries at infinity can be extended to a harmonic map 
from $X$ to $Y$ with bounded energy density if and only if the inverse 
image under the given map between the boundaries of any measure zero set
has measure zero.
\end{conj}

Note that condition on the boundary map given in the conjecture is 
exactly the assumption required for the Douady--Earle extension to 
exist.

\begin{rem} The estimates used in this paper were
strongly motivated by work of Coifman--Semmes \cite{CoiSem93}.
Although Coifman and Semmes solve a nonlinear Dirichlet problem on a bounded
domain with boundary, the norms they use make it clear that the
natural setting for their argument is on the hyperbolic
ball. On the hyperbolic ball, one no longer needs to
shrink balls as one approaches the boundary, and uniform estimates
are even easier to obtain than in the case of the bounded domain.
\end{rem}

\section{Definitions related to a map between Riemannian
manifolds}\label{defines}

We shall be considering smooth maps $u: X \to Y$. The differential of
$u$ defines a bundle map $\d u: T_*X \to T_*Y$. The {\em energy 
density} of the map $u$ is defined to be
$$
e[u] = |\d u|^2
$$

Given a bundle $F$ over $Y$ with a connection $\del$, we can pull back
both the bundle $F$ and the connection $\del$ using the map $u$,
yielding a bundle $u^*F$ over $X$ with connection $u^*\del$, which we
shall usually just denote by $\del$, for simplicity.

For example, if $u$ maps all of $X$ to a single point $y \in Y$, then
$u^*F$ is the trivial bundle $X \times F_y$ and $u^*\del$ is the
trivial connection.

The Hessian of $u$ is defined to be a bundle map $\del\d u: S^2T_*X
\to T_*Y$, where given vector fields $v$ and $w$ on $X$,
$$
\del_v\d_w u = \del_v(\d_w u) - \d_{\del_vw}u
$$
The Laplacian is then defined to be
$$
\Delta u = g^{ij}\del_i\d_ju
$$
A map $u: X \to Y$ is defined to be {\em harmonic}, if it satisfies
$$
\Delta u = 0
$$

\begin{define} A map $u: X \to Y$ is {\em strictly nondegenerate}, if
there exists a constant $\tau > 0$ such that for any $x \in X$ and $v
\in T_xX$,
\begin{equation}\label{nondeg}
\tau[u]^2  = \inf_{\substack{x \in X\\v \in T_xX, v\ne 0}} |\d u(x)|^2 -
|v|^{-2}|v\cdot\d u(x)|^2  > 0
\end{equation}
We shall call $\tau[u]$ the {\em nondegeneracy coefficient} of $u$.
\end{define}

\begin{define}
Two maps $u_0, u_1: X \to Y$ are {\em within finite distance of each other}, if
$$
\sup_{x\in X} d(u_0(x),u_1(x)) < \infty
$$
Given $1 \le p \le \infty$, the maps are 
{\em within finite $L^p$ distance of each other}, 
if the $L^p$ norm of $d(u_0(x),u_1(x))$ is finite.
\end{define}

\section{Uniqueness of homotopic harmonic maps}
Given $C^2$ homotopic maps $u_0, u_1: X \to Y$, we begin by rederiving a
well--known formula for
$\Delta d(u_0(x),u_1(x))^2$
(also, see \cite{Sam78,SchYau76}).

Given smooth homotopic maps $u_0, u_1: X \to Y$, there exits a unique section
$v$ of $u_0^*T_*Y$ satisfying
$$
u_1 = \exp_{u_0(x)} v(x)
$$
Given $0 \le s \le 1$, let
$$
u_s(x) = \exp_{u_0(x)} sv(x)
$$
Let $v_s$ be the section of $u_s^*T_*Y$ given by
$$
v_s(x) = \frac{\d}{\d s}u_s(x)
$$
Observe that $v_s(x)$, $0 \le s \le 1$, is simply the parallel
translation of $v(x)$ along the geodesic $\gamma(s) = \exp_{u_0(x)} sv(x)$.

Given $x_0 \in X$ and a unit tangent vector $e \in T_{x_0}X$,
let $x(t) \in X$, $-\delta < t < \delta$, be the
unit speed geodesic such that $x(0) = x_0$ and $x'(0) = e$.
Let
$$
\Gamma(s,t) = u_s(x(t)) = \exp_{u_0(x(t))} sv(u_0(x(t))) 
$$
and denote
\begin{align*}
\sigma &= \frac{\d}{\d s}\\
\tau &= \frac{\d}{\d t}\\
S &= \Gamma_*\sigma\\
J &= \Gamma_*\tau
\end{align*}
Let $T(s,t)$ satisfy $T(0,t) = J(0,t)$, $\del_\sigma T = 0$.
Note that
$$
S(0,t) = v(u(x(t)))
$$
Observe that
$$
0 = [\Gamma_*\sigma,\Gamma_*\tau] = \del_\sigma J - \del_\tau S
$$
Moreover, for each $t$, the curve $\Gamma(\cdot,t)$ is a constant
speed geodesic, and therefore $J$ is a Jacobi field along the geodesic
$\Gamma(\cdot,t)$, satisfying 
\begin{align*}
J(0,t) &= \Gamma_*\tau = (u_0)_*x'(t)\\
\del_\sigma J(0,t) &= \del_\sigma J(0,t) = \del_{\tau}v
\end{align*}
and the Jacobi equation 
\begin{align*}
\del_\sigma\del_\sigma J &= R(S,J)S
\end{align*}
where $R$ is the Riemann curvature tensor on $Y$.
Finally, note
\begin{align*}
\del^2_{ee}u_s &= \del_\tau\d_\tau u_s - 
\d_{\del_\tau\tau}u_s\\
&= \del_\tau J(s,0)
\end{align*}

Observe that $d(u_0(x),u_1(x))^2 = |v|^2$, so that
$$
\Delta d(u_0(x),u_1(x))^2 = 2(v\cdot\Delta v + |\del v|^2)
$$

Now consider
\begin{align*}
v\cdot\del^2_{ee}u_s &= S\cdot\del_\tau J\\
&= \frac{\d}{\d t}(S\cdot J) - J\cdot\del_\tau S
\end{align*}

Since $J$ is a Jacobi field,
\begin{align*}
\frac{\d^2}{\d s^2}(S\cdot J) &= S\cdot\del_\sigma\del_\sigma J\\
&= R(S,J)S\cdot S\\
&= 0
\end{align*}
Therefore,
$$
S\cdot J(s,t) = (S\cdot J)(0,t) + s(S\cdot\del_\sigma J)(0,t)
$$
So
\begin{align*}
S\cdot\del_\tau J(s,t) &= S\cdot \del_\tau J(0,t) - J\cdot\del_\tau S(s,t) +
J\cdot\del_\tau S(0,t)\\
&\quad+ s[\del_\tau S\cdot\del_\sigma J(0,t) + S\cdot\del_\tau\del_\sigma
J(0,t)]
\end{align*}
On other hand,
\begin{align*}
\frac{\d}{\d s}(J\cdot\del_\tau S) &= J\cdot\del_\sigma\del_\sigma J
+ \del_\sigma J\cdot\del_\sigma J\\
&= R(S,J)S\cdot J + |\del_\sigma J|^2
\end{align*}
It follows that
$$
J\cdot\del_\tau S(s,t) - J\cdot\del_\tau S(0,t) =
\int_0^s R(S,J)S\cdot J + |\del_\tau S|^2\, dr
$$
and that
\begin{align*}
v_s\cdot\del_eu_s &= S\cdot \del_\tau J\\
 &= v\cdot\del^2_{ee}u - \int_0^s
R(v_r,\d_eu_r)v_s\cdot \d_eu_r + |\del_e v_r|^2\,dr
 + s[ |\del_ev|^2 + v\cdot\del^2_{ee}v]
\end{align*}
Summing over $e$ ranging over an orthonormal frame of tangent vectors,
we get
\begin{equation*}
v_s\cdot\Delta u_s = v\cdot\Delta u + \frac{1}{2}s\Delta |v|^2
- \int_0^s R(v_r,\d^iu_r)v_r\cdot\d_iu_r + |\del v_r|^2 \,dr
\end{equation*}
and
\begin{equation}\label{distsq}
\Delta d(u_0(x),u_1(x))^2 = v_1\cdot\Delta u_1 - v\cdot\Delta u_0
+ \int_0^1 R(v_s,\d^iu_s)v_s\cdot\d_iu_s + |\del v_s|^2\,ds
\end{equation}

\begin{theo}(Theorem~\ref{unique-2})\label{uniq-2}
Assume that $\lambda_0(X) > 0$.
Then any two homotopic harmonic maps $u_0, u_1: X \to Y$ that are
within finite $L^p$ distance of each other, for some $1 < p < \infty$,
are identical.
\end{theo}

\begin{proof}
Let $\delta(x) = d(u_0(x),u_1(x))^2$. Then $\delta$ is bounded in $L^{p/2}$
by assumption and $\Delta \delta \ge 0$, by Lemma~\ref{distsq}. It
follows that $\del(\delta^{p/4})$ is in $L^2$ and that
\begin{align*}
\int \delta^{p/2} &\le \lambda_0(X)^{-1}\int |\del\delta^{p/4}|^2\\
&= C\int \delta^{(p-2)/2}(-\Delta\delta)\\
&\le 0
\end{align*}
\end{proof}

\begin{rem} Theorem~\ref{uniq-2} and its proof hold, if $Y$ has nonpositive
sectional curvature.
\end{rem}

\begin{theo}(Theorem~\ref{unique-1})
Let $u_0: X \to Y$ be a $C^2$ strictly nondegenerate harmonic map. Any
$C^2$ harmonic map $u_1: X \to Y$ that is homotopic to
and within finite distance of $u_0$ is equal to $u_0$.
\end{theo}

\begin{proof} Given  homotopic maps $u_0$ and $u_1$, there
exists a $C^2$ section $v$ of $u_0^*T_*Y$ such that $u_1 = \exp_{u_0}
v$.

The function $\delta(x) = d(u_0(x),u_1(x))^2$ is $C^2$.
 By the generalized maximum principle, Lemma~\ref{max}, there exists a sequence $x_i \in X$ such that
\begin{align*}
\lim_i \delta(x_i) &= \sup_{x \in X} \delta(x)\\
\limsup_i \Delta \delta(x_i) &\le 0
\end{align*}
On the other hand, by \eqref{distsq},  Lemma~\ref{indexest}, and the
strict nondegeneracy of $u_0$,
\begin{align*}
\Delta\delta &\ge \kappa^{-1} |v|(\tanh \kappa^{-1}|v|)[|\d u_0|^2 -
|v|^{-2}|v\cdot\d u_0|^2]\\
&\ge \kappa \tau^2 \sqrt{\delta}\tanh \kappa \sqrt{\delta}
\end{align*}
Therefore,
\begin{align*}
\sup_{x\in X} \kappa^{-1} \tau^2 \sqrt{\delta}\tanh \kappa^{-1} \sqrt{\delta}
&= \lim_i \kappa^{-1} \tau^2 \sqrt{\delta(x_i)}\tanh \kappa^{-1}
\sqrt{\delta(x_i)}\\
&\le \limsup_i \Delta\delta(x_i)\\
&\le 0
\end{align*}
We conclude that $\delta$ vanishes identically.
\end{proof}

\section{A flow for maps between Riemannian manifolds}

We shall demonstrate the existence of a harmonic map by constructing
a flow of maps that converges to a harmonic map.
It is probably
possible to use the harmonic map heat flow, but we use
a different flow that appears particularly well--suited to maps 
that have infinite total energy but bounded energy density.

The flow we use is the following: Given a $1$--parameter family of
maps $u(t): X \to Y$, $0 \le t < \infty$, let $v = \d_tu$ denote the
velocity field. Note that $v(t)$ is a section of $u(t)^*T_*Y$. The
equation for the flow is
$$
\Delta v - R(\d^iu,v)\d_iu = -\Delta u
$$
where $R$ is the Riemann curvature of $Y$ and $\Delta$ denotes the
Laplacian(s) defined using the naturally induced connections on the
corresponding bundles.
The choice of this flow is motivated by the observation that
$$
\Delta \d_{t}u = \del_{t}\Delta u
$$
Therefore, if $u(t)$ is a solution for the flow, it satisfies
$$
\del_{t}\Delta u(t) = -\Delta u(t)
$$

Under the assumptions we impose, the linear differential operator
$\Delta - R(\d^iu,\cdot)\d_iu$ is always invertible and has a
bounded inverse. This leads to the existence and uniqueness of the
flow for small time.

On other hand, using the generalized maximum principle for complete
Riemannian manifolds and local $L^{p}$ bounds for
elliptic PDE's, it is easy to obtain {\it a priori} bounds for the
flow that imply long time existence and convergence to a harmonic map as
$t \to \infty$.

\begin{rem} The flow defined above was introduced in the doctoral
dissertation of J. P. Anderson \cite{and91}, who used it to
deform a
quasiisometry between closed Riemannian manifolds into a harmonic
diffeomorphism.
\end{rem}

\section{Basic analytic definitions and estimates}
\label{basic}

The distance between two points $x$ and $y$  on a Riemannian 
manifold $X$
will always be denoted $d(x,y)$. A geodesic ball of radius $r$ 
centered at $p$ will be denoted $B(p,r)$. Let $\widetilde{B}(p,r)$ 
denote the universal cover of $B(p,r)$ and $\Bhat(p,r) \subset 
\widetilde{B}(p,r)$ be the ball of radius $r$ centered at one of the 
inverse images of $p$, using the pullback Riemannian metric.
Let $\phi_{p}: \Bhat(p,r) \to B(p,r)$ denote the canonical projection 
map.

Recall the Bishop--Gromov volume comparison theorem (see 
\cite{Cha84}):
\begin{lem}
Given a complete Riemannian manifold $X$ with Ricci curvature bounded 
from below by $-(n-1)\kappa^{2}$, $p \in X$, and $0 < r < R$,
$$
\frac{\vol(\Bhat(p, R))}{\vol(\Bhat(p,r))} \le 
\frac{v_{-\kappa^{2}}(R)}{v_{-\kappa^{2}}(r)}
$$
where $v_{-\kappa^{2}}(r)$ is the volume of a geodesic ball with 
radius $r$ in the space of constant sectional curvature 
$-\kappa^{2}$. 
\end{lem}

\begin{define}
Given a complete Riemannian manifold $X$ with Ricci curvature bounded 
from below by $-(n-1)\kappa^{2} \le 0$ and $x \in X$, we define the 
{\em Sobolev radius $\rho_{S}(x)$} to be the largest $r > 0$ such that
$$
\vol(\Bhat(x,r)) \ge \frac{1}{2}v_{-\kappa^{2}}(r)
$$
This definition is motivated by the theorem of Gromov (see 
\cite{Cha84}), which gives a lower bound on the isoperimetric constant
and therefore an upper bound on the Sobolev constant on the domain 
$\Bhat(x,\rho_{S}(x))$.
\end{define}

Recall that we always assume that the Ricci curvature of $X$ is 
bounded from below by $-(n-1)\kappa^{2}$ and that the Sobolev radius 
is always bounded from below by $2\rho > 0$.
As a consequence, we obtain the following key lemmas:
\begin{smallcor} There exists $N(n,\kappa) > 0$ such that given any $x \in
X$, there exist $x_1, \dots, x_N \in X$ such that
$$
\widehat{B}(x,2\rho) \subset \widehat{B}(x_1,\rho) \cup \dots \cup
\widehat{B}(x_N,\rho)
$$
\end{smallcor}

Let $E$ be a vector bundle over $X$ with an inner product.
Given a section $v$ of $E$, $1 \le p < \infty$, 
we shall use the following ``local'' $L^p$ norms:
\begin{align*} 
\|v\|_{\phat} &= 
\sup_{x \in X}\left(\int_{\widehat{B}(x,\rho)} 
|v\circ\phi_x|^p\right)^{1/p}
\end{align*}
As mentioned before, the idea of using norms like this is motivated by a
similar approach used by Coifman--Semmes \cite{CoiSem93} to solve a Dirichlet
problem for a nonlinear elliptic PDE on a bounded domain.

\begin{smallcor} \label{integrals}
Given any locally integrable function $f: X \to R$ and $x \in
X$,
$$
\int_{\widehat{B}(x,2\rho)} f\circ\phi_x \le N\sup_{y \in X}
\int_{\widehat{B} (y,\rho)} f\circ\phi_x 
$$
In particular, for any locally $L^p$ function $f$ and $x \in X$,
\begin{align*}
\left(\int_{\Bhat(x,2\rho)} |f|^p\right)^{1/p} \le N^{1/p}\|f\|_{\phat}
\end{align*}
\end{smallcor}

The key analytic estimate we will use is obtained using a classical 
technique known as Moser iteration. The version we use is taken from 
Appendix C of \cite{yan92c}.

\begin{lem}\label{moser}
Let $E$ be a vector bundle over $X$ with an inner product 
and a compatible connection. Let $R: E \to E$ be a 
symmetric linear bundle map and
$$
K_{-}(x) = \inf_{v \in E_{x}\backslash\{0\}} \frac{v\cdot R(x)v}{|v|^{2}}
$$
Given $2 \le p, q \le \infty$ and $n/2 < s 
\le \infty$, there exists
$C > 0$ such that given sections $v$ of $E$ and $f$ of 
$T_{*}X\otimes E$ that satisfy
\begin{align*}
\|v\|_{\phat} &< \infty\\
\|f\|_{\qhat} &< \infty\\
\|K_{-}\|_{\infty} &< \infty\\
-\del^{i}\del_{i}v + Rv &= \del\cdot f
\end{align*}
then if $p < q < n$,
\begin{align}
\|v\|_{\hat{r}} &\le C[(1 + 
\|K_{-}\|_{\infty})^{\frac{n}{2}\left(\frac{1}{p}-\frac{1}{r}\right)}
\|v\|_{\phat} + \|f\|_{\hat{q}}]
\end{align}
where
$$
r = \frac{qn}{n-q}
$$
and if $q > n$,
\begin{align}
\|v\|_{\infty} &\le C[(1 + 
\|K_{-}\|_{\infty})^{\frac{n}{2p}}\|v\|_{\phat} + 
\|f\|_{\hat{q}}]
\end{align}
\end{lem}

\section{A priori estimates for the harmonic map flow}

\begin{lem} \label{p-estime} Assume $\lambda_{0}(X) > 0$.
 Given $p \ge 2$, there exists $C > 0$
such
that the following holds: Given a smooth map $u: X \to Y$
such that $\|\Delta u\|_{p} < \infty$
and an $L^{p}$--bounded section $w$ of $u^*T_*Y$, there exists a
unique $L^{p}$--bounded section $v$ of $u^*T_*Y$ satisfying
\begin{equation}\label{linearized}
-\Delta v + R(\d^iu,v)\d_iu = w
\end{equation}
Moreover, $v$ satisfies the following estimates:
\begin{align*}
\|v\|_p &\le C\|w\|_p\\
\|\del v\|_{\qhat} &\le C\|w\|_{p}\\
\end{align*}
where
\begin{align*}
q &= \begin{cases}
\frac{pn}{n-p} & p < n\\
\infty & p > n
\end{cases}
\end{align*}
\end{lem}

\begin{proof}
First, note that
$$
v\cdot w = v\cdot(-\Delta v) + v\cdot R(\d^iu,v)\d_iu \ge v\cdot(-\Delta v)
$$
Therefore,
\begin{align*}
\int |v|^p &\le \lambda_0^{-1}\int |\del(|v|^{p/2})|^2\\
&\le \frac{p^2}{4(p-1)\lambda_0}\int |v|^{p-2}|\del v|^2 +
(p-2)|v|^{p-4}|v\cdot\del v|^2\\
&= C\int |v|^{p-2}v\cdot(-\Delta v)\\
&\le C\int |v|^{p-2}v\cdot w\\
&\le C\|v\|_p^{p-1}\|w\|_p
\end{align*}
This proves the first inequality and the uniqueness of $v$.

Next, observe that using the equation \eqref{linearized}, we get
$$
\|\del v\|_{\hat{2}} \le C\|v\|_{\hat{2}} \le C\|v\|_{\phat} \le 
C\|w\|_{\phat}
$$
Differentiating and manipulating \eqref{linearized}, we obtain
\begin{multline}\label{linearized'}
-\del^{j}\del_{j}\del_{i}v^{\alpha}
- R^{\alpha}{}_{\beta\gamma\delta}
\d_{i}u^{\gamma}\d_{j}u^{\delta}\del^{j}v^{\beta}
+ \widehat{R}_{i}{}^{j}\del_{j}v^{\alpha}\\
= 
- \del^{j}(R^{\alpha}{}_{\beta\gamma\delta}\d_{j}u^{\gamma}\d_{i}u^{\delta}v^{\beta})
+ \del_{i}(R^{\alpha}{}_{\beta\gamma\delta}\d^{p}u^{\beta}\d_{p}u^{\delta}v^{\gamma}
+ w)
\end{multline}
The bound for $\del v$ now follows from the local $L^{2}$ bound on 
$\del v$ and Lemma~\ref{moser}. When applying Lemma~\ref{moser}, note 
that we can throw away the term containing the curvature of $Y$, 
because it is positive definite.
 
To obtain existence, consider the equation on a sequence of domains
with smooth boundary in $X$ that exhausts $X$. On each domain the
positivity of the potential implies the
existence and uniqueness of a solution with Dirichlet boundary
conditions, i.e. boundary value equal to
zero. The estimates imply that on any compact domain this sequence of
solutions has a subsequence that converges uniformly to a solution.
 Using a diagonal argument, we obtain a subsequence that converges
 on any compact subset of $X$ to a global solution of the equation.
\end{proof}

Integrating the estimates above yields
\begin{lem} \label{p-apriori} Assume $\lambda_{0}(X) > 0$. Given $1 < p < \infty$, let
  $C > 0$ and $q > n$ be as given by Lemma~\ref{estime}.
Let $u(t): X \to Y$, $0 \le t < T$ be a $1$--parameter family of maps
satisfying
\begin{align*}
\|v(t)\|_p &< \infty\\
\|\Delta u(t)\|_{p} &< \infty\\
-\Delta v(t) + R(\d^iu(t),v(t))\d_iu(t) &= \Delta u(t)
\end{align*}
where $v(t) = \d_tu(t)$.
Then the following estimates hold for any $0 \le t < T$:
\begin{align*}
\|\Delta u(t)\|_{p} &= \|\Delta u(0)\|_{p} e^{-t}\\
\|d(u(t), u(0))\|_{p} &\le C\|\Delta u(0)\|_{p}(1 - e^{-t})\\
\|\d u(t)\|_{\qhat} &\le \|\d u(0)\|_{\qhat} + C\|\Delta 
u(0)\|_{p}(1-e^{-t})
\end{align*}
\end{lem}

With an $L^{\infty}$ bound on $\Delta u(t)$, it's necessary to assume 
a lower bound on the nondegeneracy of $u(t)$.

\begin{lem}\label{estime} There exists $C > 0$
such that the following holds:
Given a map $u: X \to Y$ with $L^{\infty}$--bounded $\d u$ and
nondegeneracy coefficient $\tau = \tau[u]$ and a bounded
smooth section $w$ of
$u^*T_*Y$, there exists a  unique bounded smooth section $v$ of $u^*T_*Y$ satisfying
\begin{equation}
-\Delta v + R(\d^iu,v)\d_iu = w
\end{equation}
Moreover, $v$ satisfies the following estimates:
\begin{align*}
\|v\|_\infty &\le \left(\frac{\kappa}{\tau}\right)^{2}\|w\|_\infty\\
\|\del v\|_{\infty} &\le C\|w\|_{\infty}
\end{align*}
where the sectional curvature of $Y$ is bounded between $-\kappa^{2}$ 
and $\kappa^{-2}$.
\end{lem}

\begin{proof}
Uniqueness of $v$ follows by applying the generalized maximum
principle, Lemma~\ref{max}, to the elliptic inequality satisfied by
$|v|^2$. The maximum principle also yields the bound on
$\|v\|_\infty$. The $L^{\infty}$ bound on $\del v$ follows 
as in the proof of Lemma~\ref{p-estime} by applying Lemma~\ref{moser}
to \eqref{linearized'}.

The existence of a solution $v$ follows by the same argument as in the 
proof of Lemma~\ref{p-estime}.
\end{proof}

Integrating the estimates in Lemma~\ref{estime} yields
\begin{lem} \label{apriori} 
Let  $C > 0$ be as given by Lemma~\ref{estime}.
Given $\tau > 0$, let $u(t): X \to Y$, $0 \le t < T$ be a $1$--parameter family of maps
satisfying
\begin{align*}
\tau[u(t)] &\ge \tau\\
\|v(t)\|_\infty &< \infty\\
\|\Delta u(t)\|_{\infty} &< \infty\\
-\Delta v(t) + R(\d^iu(t),v(t))\d_iu(t) &= \Delta u(t)
\end{align*}
where $v(t) = \d_tu(t)$.
Then the following estimates hold for any $0 \le t < T$:
\begin{align*}
|\Delta u(t)| &= |\Delta u(0)| e^{-t}\\
\|\d u(t)\|_{\infty}
&\le \|\d u(0)\|_{\infty} + C\|\Delta u(0)\|_{\infty}(1-e^{-t})
\end{align*}
\end{lem}

\section{Existence theorems for the harmonic map flow}

We begin with existence of the flow for small time:

\begin{theo} \label{flow}Let $u_0: X \to Y$ be a map satisfying
$$
\|\d u_0\|_{\infty}, \|\Delta u_0\|_\infty < \infty
$$
and
$$
\tau[u_0] > 0
$$
Then there exists $T > 0$ and a unique $1$--parameter family of maps
$u: X \times [0,T) \to Y$ and $v \in \Gamma(u^*T_*Y)$ satisfying 
the harmonic map flow:
\begin{align*}
\del_t u &= v\\
\|v\|_{\infty} &< \infty\\
-\Delta v + R(\d^iu,v)\d_iu&=\Delta u\\
u(0) &= u_0
\end{align*}
Moreover, if
\begin{equation}\label{long}
\limsup_{t\to T} \tau[u(t)] > 0
\end{equation}
then the flow extends beyond $t = T$.
\end{theo}

\begin{proof}
The existence and uniqueness of a solution for sufficiently small $T
> 0$ follows from a straightforward application of the usual existence
and uniqueness theorem for an ODE on a Banach manifold. We sketch
the argument.

Given $C > 0$, the space of maps $u: X \to Y$ satisfying
\begin{align*}
\tau[u] &> C^{-1}\\
\|d(u_0,u)\|_\infty &< C\\
\|\d u\|_{\infty} &< C\\
\|\Delta u\|_\infty &< C
\end{align*}
is an open Banach manifold $\mathcal{B}$. We leave to the reader the
exercise of constructing coordinate charts for $\mathcal{B}$ and
verifying that given $u \in \mathcal{B}$, the tangent space
$T_u\mathcal{B}$ consists of sections $v$ of $u^*T_*Y$ satisfying
$$
\|v\|_\infty, \|\del v\|_{\infty}, \|\Delta
v\|_\infty < \infty
$$
Let $F$ be the section of $T_*\mathcal{B}$ where $v = F(u)$
is the unique section of $u^*T_*Y$ satisfying
$$
\Delta v - R(\d^iu,v)\d_iu = -\Delta u
$$
By Lemma~\ref{estime}, $F$ is a well-defined differentiable section
of $T_*\mathcal{B}$. By the existence and uniqueness theorem for
ODE's on a Banach manifold \cite{Lan95},
there exists for sufficiently small $T > 0$ a unique solution $u: [0,T) \to \mathcal{B}$ of
$$
\frac{d u}{d t} = F(u)
$$
If there exists $\tau > 0$ such that
$$
\tau[u(t)] \ge \tau
$$
for all $0 \le t < T$,
 then the estimates obtained in Lemmas~\ref{estime} and 
\ref{apriori} imply that $u(t) \to u(T) \in \mathcal{B}$, for $C$ 
suitably large. Then the local existence theorem implies that the
solution can be extended beyond $t = T$. \end{proof}

\begin{theo}(Theorem~\ref{exist-1})\label{har}
Given manifold $X$, $Y$, and  $\tau > 0$, there exists $\e > 0$,
such that if $u_0: X \to Y$ has nondegeneracy coefficient $\tau[u_0]
\ge \tau$
and $$
\|\Delta u_0\|_\infty < \e
$$
then there exists a unique homotopic harmonic map $u_1: X \to Y$
within finite distance of $u_0$.
\end{theo}

\begin{proof} By Theorem~\ref{flow}, there exists a harmonic map flow 
$u(t): X \to Y$ with $u(0) = u_{0}$. By Lemma~\ref{apriori}, if $\e > 
0$ is chosen sufficiently small, then $\tau[u(t)] \ge \tau/2$ for any 
$t$. Therefore, the flow must exist for all $t > 0$.
 The bounds given in Lemma~\ref{estime} imply that as $t
\to \infty$, the map $u(t)$ converges to a limiting map $u_\infty$.
Moreover, the limiting map must be harmonic and bounded distance
from $u_0$.
\end{proof}

In the same way, Lemma~\ref{p-apriori} leads to the following
\begin{theo} (Theorem~\ref{exist-2})
Assume that
$$
\lambda_0(X) > 0
$$
Given $2 \le p < \infty$, let
$$
q = \begin{cases}
\frac{pn}{n-p} & 1 < p < n\\
\infty & p > n
\end{cases}
$$
Let $u_0: X \to Y$ be a smooth map satisfying
\begin{align*}
\|\d u_0\|_{\qhat} &< \infty\\
\|\Delta u_0\|_p &< \infty\\
\end{align*}
Then there a unique $1$--parameter family of maps
$u: X \times [0,\infty) \to Y$ and $v \in \Gamma(u^*T_*Y)$
satisfying:
\begin{align*}
\del_t u &= v\\
\|v\|_p &< \infty\\
-\Delta v + R(\d^iu,v)\d_iu&=\Delta u\\
u(0) &= u_0
\end{align*}
Moreover, there exists a smooth limit
$$
u_\infty = \lim_{t\to\infty} u(t)
$$
that is a harmonic map and is within finite $L^p$ distance of the
map $u_0$. \end{theo}

What is not obvious is that given two Hadamard manifolds $X$ and $Y$,
whether there exists any nontrivial map $u_0$ that has $\Delta u_0$
bounded in $L^p$ for some $p$. This is, in fact, unlikely to occur in
the general situation. Understanding when this is possible would
extend the results of Li--Tam to Hadamard manifolds.

\appendix

\section{A generalized maximum principle for complete
Riemannian manifolds}

\begin{lem}\label{cutoff} Given $\kappa, r > 0$, there exists $C > 0$ such that
for any complete Riemannian manifold with Ricci curvature bounded from 
below by $-(n-1)\kappa^{2}$, $p \in X$, there exists a compactly 
supported function $\chi$ on $\Bhat(p, 2r)$ such that the following 
hold:
\begin{align*}
\chi(x) &= 1,\ x \in \Bhat(p,r)\\
\|\del\chi\|_{\infty} &\le C\\
\Delta\chi(x) &\le C
\end{align*}
\end{lem}

\begin{proof} A version of the Bishop--Gromov comparison inequality 
states that the function $\rho(x) = d(p, x)$ satisfies
$$
\Delta \rho \le (n-1)\kappa \coth \kappa \rho
$$
in a generalized sense. This is easily proven whenever there is a 
unique minimal geodesic joining $p$ to $x$. If there is not, then one 
simply observes that in a neighborhood of $x$, $\rho$ can be 
represented as the minimum of a finite number of smooth functions 
that satisfy the inequality.

It suffices to set $\chi(x) = \psi(\rho(x))$, where $\psi$ is a 
suitably chosen smooth compactly supported function on the real line.
\end{proof}

\begin{lem}(Omori--Yau maximum principle, \cite{Omo67,Yau75},
 also see \cite{Aub82})
\label{max}
Given a complete Riemannian manifold $X$ with Ricci curvature bounded
from below
and  any $C^2$ function $f: X \to \R$ that is bounded from above,
there exists a sequence of points $x_k \in X$ such that
\begin{align*}
\lim_k f(x_k) &= \sup f\\
\lim_k |\del f(x_k)| &= 0\\
\limsup_k \Delta f(x_k) &\le 0
\end{align*}
\end{lem}

\begin{proof} It suffices to prove that there exists a constant $C
> 0$ such that 
given any $C^2$ function $f$ that is bounded from above and $\e > 0$,
there exists $x \in X$ such that
\begin{equation}\label{ebds}
\begin{split}
f(x) &\ge \sup f - \e\\
|\del f(x)| &\le C\e\\
\Delta f(x) &\le C\e
\end{split}
\end{equation}
Given $\e > 0$ and a $C^2$ function $f$ that is bounded from above,
there exists $x' \in X$ satisfying $f(x') \ge \sup f - \e$.
Consider the function $\hat{f}: \Bhat(x',2) \to \R$, where
$$
\hat{f}(z) = f(\phi_{x'}(z)) - \e[1-\chi(z)]
$$
where $\chi$ is the function given by Lemma~\ref{cutoff}.
Note that $\hat{f}$ must achieve a maximum at some $z' \in B_0$.
Therefore,
\begin{align*}
\del\hat{f}(z') &= 0\\
\Delta\hat{f}(z') &\le 0
\end{align*}
Therefore, \eqref{ebds} hold if
$x = \phi_{x'}(z')$.
\end{proof}

\section{Lower bounds for the index of a geodesic}

Let $Y$ be a complete Riemannian manifold with strictly negative curvature.
Let $\gamma: [0,1] \to Y$ be a constant speed geodesic. Let $T$ be a
nonzero tangent vector at $\gamma(0)$.
We want to minimize the quantity
$$
E[J] = \int_0^1 R(\gamma',J)\gamma'\cdot J + |\del_{\gamma'}J|^2\,ds
$$
over all Jacobi fields $J$ along $\gamma$ such that $J(0) = T$. Since
the functional $E$ is a positive definite quadratic function on a finite
dimensional affine space, there is a unique minimum.

 If $J$ is the
minimizing Jacobi field, then given any Jacobi field $\dot{J}$ along
$\gamma$ such that $\dot{J}(0) = 0$,
\begin{align*}
0 &= E'[J]\dot{J}\\
&= 2\int_0^1 R(\gamma',J)\gamma'\cdot \dot{J} + J'\cdot\dot{J}'\,ds\\
&= \left. 2\dot{J}\cdot J'\right|_0^1\\
&= 2\dot{J}\cdot J'(1)
\end{align*}
Therefore, $J'(1) = 0$ and
\begin{align*}
\inf_{\substack{
\text{$J$ Jacobi}\\
J(0) = T
}}
\int_0^1 R(\gamma',J)\gamma'\cdot J + |\del_{\gamma'}J|^2\,ds
&=
\inf_{\substack{
\text{$J$ Jacobi}\\
J(0) = T, J'(1) = 0
}}
\int_0^1 R(\gamma',J)\gamma'\cdot J + |\del_{\gamma'}J|^2\,ds\\
&=
\inf_{\substack{
J(0) = T\\ J'(1) = 0}}
\int_0^1 R(\gamma',J)\gamma'\cdot J + |\del_{\gamma'}J|^2\,ds\\
&\ge
\inf_{\substack{
J(0) = T\\ J'(1) = 0}}
\int_0^1 \alpha^2[|\gamma'|^2|J|^2-(\gamma'\cdot J)^2] + |J'|^2\,ds
\end{align*}
where the sectional curvature of $Y$ is bounded from above by
$-\kappa^{-2}$.

We can, however, compute the last infimum explicitly, since it is
achieved by the vector field
$$
J(s) = |\gamma'|^{-2}(\gamma'\cdot T)\gamma'
+ (T-|\gamma'|^2(\gamma'\cdot T)\gamma')(\cosh s\kappa^{-1} |\gamma'| - 
\tanh \kappa^{-1}|\gamma'|\sinh s\kappa^{-1}|\gamma'|)
$$
Differentiating this, we get
$$
J'(s) = \kappa^{-1}|\gamma'|
(T-|\gamma'|^{-2}(\gamma'\cdot T)\gamma')(\sinh s\kappa^{-1} |\gamma'| - 
\tanh \kappa^{-1}|\gamma'|\cosh s\kappa^{-1}|\gamma'|)
$$
Then
\begin{align*}
\int_0^1 \kappa^{-2}|(\gamma'|^{-2}|J|^2 - (\gamma'\cdot J)^2)
+ |J'|^2\,ds
&= \left.J\cdot J'\right|_0^1\\
&= \kappa^{-1}|\gamma'|(\tanh \kappa^{-1}|\gamma'|)
(|T|^2-|\gamma'|^{-2}(\gamma'\cdot T)^2)
\end{align*}

We conclude that
\begin{lem}\label{indexest}
Given a Riemannian manifold $Y$ with sectional curvature bounded from
above by $-\kappa^{-2} < 0$, a constant speed geodesic $\gamma: [0,1] \to
Y$, and a nonzero tangent vector $T$ at $\gamma(0)$, then given any
Jacobi field $J$ such that $J(0) = T$,
$$
\int_0^1 R(\gamma',J)\gamma'\cdot J + |\del_{\gamma'}J|^2\,ds
\ge \kappa^{-1}|\gamma'|(\tanh \kappa^{-1}|\gamma'|)
[|T|^2-|\gamma'|^{-2}(\gamma'\cdot T)^2]
$$
\end{lem}

\section{The Douady--Earle extension} \label{de}

In this section let $B$ denote the unit ball in $\R^n$.
Let
$$
g = \frac{4|dx|^2}{(1-|x|^2)^2}
$$
be the standard hyperbolic metric on $B$. Throughout
this section, norms denoted $|\cdot|$ and $\|\cdot\|$ are with respect
to the Euclidean metric, and norms denoted $|\cdot|_g$ and
$\|\cdot\|_g$ are with respect to the hyperbolic metric.

Douady and Earle \cite{DouEar86} proved the following:
\begin{theo}
Given a homeomorphism $\hat{u}: \d B \to \d B$, there exists a unique
map $\mathcal{E}[\hat{u}] = u: B \to B $ such that for any $\theta \in \d B$,
$$
\lim_{x \to \theta} u(x) = \hat{u}(\theta)
$$
and given any $x \in B$,
$$
\int_{\d B} \phi_{u(x)}\circ\hat{u}\circ\phi^{-1}_x(\theta)\,d\theta =
0
$$
where for any $z \in B$, $\phi_z: \overline{B} \to \overline{B}$ is a M\"obius
transformation such that $\phi_z(z) = 0$.

The map $u$ is real analytic.
Given any M\"obius transformations $\phi, \psi$,
$$
\mathcal{E}[\phi\circ\hat{u}\circ\psi] =
\phi\circ\mathcal{E}[\hat{u}]\circ\psi
$$
\end{theo}

Douady and Earle also showed that the value of the extended map and
its derivatives at any given point depend continuously on the
boundary data. In particular,
\begin{lem}\label{de-bound} Given a homeomorphism $\hat{u}_0: \d B \to \d B$, $k > 0$,
and  any
$\e > 0$,
there exists $\delta > 0$ such that
given any homeomorphism $\hat{u}: \d B \to \d B$ satisfying
$$
\|\hat{u}-\hat{u}_0\|_\infty < \delta
$$
we have
$$
|\del^k\mathcal{E}[\hat{u}](0) - \del^k\mathcal{E}[\hat{u}_0](0)| < \e
$$
\end{lem}

\begin{lem} \label{compact}Given $K > 1$ and any sequence $\hat{u}_{i}: \d B \to \d 
B$ of $K$--quasisymmetric maps such that the extended maps $u_{i} = 
\mathcal{E}[\hat{u}_{i}]$ satisfy $u_{i}(0) = 0$ has a uniformly 
convergent subsequence.
\end{lem}

\begin{proof} Suppose not. Then there exists M\"obius 
transformations $\phi_{i}$ and a subsequence $\hat{u}_{i}$ such that 
$\phi_{i}\circ\hat{u}_{i}$ does converge uniformly. Therefore, 
the sequence $\phi_{i}\circ u_{i}(0)$ converges to a limit inside $B$.

On the other hand, no subsequence
of $\phi_{i}$ converges. Therefore, there exists a subsequence 
$\phi_{i}$ such that $\phi_{i}\circ u_{i}(0) = \phi_{i}(0) \to \d B$. 
This is a contradiction.
\end{proof}

\begin{smallcor}\label{qiext}
 Given $\e > 0$, there exists $\delta > 0$ such that if
$\hat{u}: \d B \to \d B$ is $(1+\delta)$--quasisymmetric, then $u =
\mathcal{E}[\hat{u}]$ is a $(1+\e)$--quasiisometry satisfying
\begin{align*}
\|\del\d u\|_{g,\infty} &\le \e\\
\end{align*}
\end{smallcor}

\begin{proof} Suppose not. Then for each $k \in \Z^+$, there exists a
$(1+k^{-1})$--quasisymmetric map
$\hat{u}_k: \d B \to \d B$ and $x_k \in D$ such that  either
\begin{align}\label{cond1}
|u_k^*g(x_k) - g|_g &> 2\e\\
\intertext{or}
|\del\d u_k(x_k)|_g &> \e \label{cond2}
\end{align}
By composing with M\"obius transformations (which are isometries of
the metric $g$), we may assume that for all
$k$, $x_k = u_k(x_k) = 0$. By Lemma~\ref{compact} we may restrict to a
subsequence $\hat{u_k}$ that converges uniformly to a M\"obius
transformation. Therefore, by Lemma~\ref{de-bound},
$u_k^*g(0)$ must converge to
$g$ and $\del\d u_k(0)$ to $0$. This contradicts \eqref{cond1} and \eqref{cond2}.
\end{proof}

In dimension $2$ Douady--Earle prove more. They show that if the boundary map is quasisymmetric, then the extension is always quasiconformal. In particular,
\begin{theo} 
 Given $K > 1$, there exists $K^* \ge K$
such that if $\hat{u}: S^1 \to S^1$ is $K$--quasisymmetric, then
$\mathcal{E}[\hat{u}]: D \to D$ is $K^*$--quasiconformal.
\end{theo}

\bibliographystyle{amsplain}

\providecommand{\bysame}{\leavevmode\hbox to3em{\hrulefill}\thinspace}

\end{document}